\documentclass{amsart}

\usepackage{graphicx}

\newtheorem{theorem}{Theorem}[section]
\newtheorem{lemma}[theorem]{Lemma}

\theoremstyle{definition}
\newtheorem{definition}[theorem]{Definition}

\newcommand{\ind}{\mathrm{Ind}}
\newcommand{\nhs}{\dot{N}}
\newcommand{\susp}{\mathrm{susp}}
\newcommand{\sm}{\setminus}
\newcommand{\bu}{\item[$\ast$]}

\newcommand{\hideS}[1]
  {\makebox[0cm]{\scriptsize $#1$}}
\newcommand{\hideT}[1]
  {\makebox[0cm]{\tiny $#1$}}
\newcommand{\wideEqual}
  {\makebox[1cm]{$=$}}
\newcommand{\wideCup}
  {\makebox[1.5cm]{$\cup$}}

\begin{document}

\title
  [Independence complexes of claw-free graphs]
  {Independence complexes of claw-free graphs}

\author{Alexander Engstr\"om}

\address{Department of Computer Science, 
  Eidgen\"ossische Technische
  Hochschule, Z\"urich, Switzerland.} 

\email{engstroa@inf.ethz.ch}

\thanks{Research supported by ETH and Swiss 
  National Science Foundation Grant 
  PP002-102738/1}


\date\today


\begin{abstract}
  We study the class of independence complexes of claw-free
  graphs. The main theorem give good bounds
  on the connectivity of these complexes, given
  bounds for a few subcomplexes of the same class. 
  Two applications are presented.
  Firstly, we show that the independence complex of a 
  claw-free graph with $n$ vertices and maximal degree 
  $d$ is $(cn/d+\varepsilon)$--connected,
  where $c=2/3$. This can be compared with the result of
  Szab\'o and Tardos that $c=1/2$ is optimal with
  no restrictions on the graphs. Secondly, we calculate
  the connectivity of a family of complexes used
  in Babson and Kozlov's proof of Lov\'asz conjecture.
\end{abstract}

\maketitle

\section{Introduction}
  The independence complex is a good structure
  for transferring graph coloring problems
  to combinatorial topology. Usually the
  topological statements to investigate will
  be about connectivity. In this paper we study
  the connectivity of independence complexes
  of claw-free graphs.

  First let us fix notation and
  introduce some tools.
\subsection{Graphs}
  All graphs are finite and simple. For a 
  graph $G$ the
  edge set is $E(G)$ and the vertex set $V(G)$.
  A complete graph has edges between all 
  vertices. The complement of $G$ is called 
  $\overline{G}$.
  The induced subgraph of $G$ on 
  $U\subseteq V(G)$ is denoted $G[U]$, and
  $G\sm U = G[V(G)\sm U]$. A set
  $I\subseteq V(G)$ is independent if $G[I]$ 
  lacks edges. The set of vertices of a 
  graph $G$ 
  with edges to a vertex $v$, is the 
  neighborhood of $v$. It is called 
  $N_G(v)$, or just $N(v)$. And $\nhs(v)=
  N(v)\cup\{v\}$. 
\subsection{Topological tools}\label{ss:tt}  
  All topological tools used are standard.
  For proofs and further references see 
  Bj\"orner's survey \cite{B} chapters 9--10.
  A topological space $T$ is $n$-connected if
  for all $0\leq i\leq n$ any map from the
  $i$--sphere to $T$ can be extended 
  to a map from the $(i+1)$--ball to $T$. 
  Arcwise connected and $0$--connected is the 
  same. Define all non-empty spaces to be
  $(-1)$-connected, and all spaces to be 
  $n$--connected for $n\leq -2$.  
  These lemmas will be used several times:
  \begin{lemma}[Corollary of Theorem 10.6 
  \cite{B},Theorem 1.1 \cite{BLVZ}]
  \label{lemma:cf1}
  If $\Delta_1, \Delta_2, \ldots
  \Delta_k$ are $n$--connected simplicial
  complexes and 
  $\cap_{i\in I}\Delta_i$
  is  $(n-1)$--connected
  for any $\emptyset\neq I \subseteq \{1,2,\dots
  k\}$ then $\cup_{i=1}^{k}\Delta_i$ is
  $n$--connected.
  \end{lemma}
  \begin{lemma}[Theorem 10.4 \cite{B}]
  \label{lemma:cf2}
  If $\Delta_0, \Delta_1, \ldots
  \Delta_k$ are contractible simplicial complexes
  and $\Delta_i\cap\Delta_j\subseteq \Delta_0$
  for all $1\leq i<j\leq k$ then
  $\cup_{i=0}^k\Delta_i=\vee_{i=1}^k 
  \susp (\Delta_0\cap\Delta_i).$
  \end{lemma}
  If $\Delta$ is a simplicial complex with
  vertex set $V$ and $U\subseteq V$, then the 
  induced
  subcomplex is $\Delta[U]=\{\sigma\in\Delta
  \mid \sigma\subseteq U\}$.
\section{Independence complexes of claw-free
         graphs}
\subsection{Claw-free graphs}
  A \emph{claw} is four vertices $u,v_1,v_2,v_3$
  with edges from $u$ to $v_1,v_2,v_3$, but
  no edges among $v_1,v_2,v_3$. A graph is
  \emph{claw-free} if there are no induced
  subgraphs which are claws. An equivalent
  definition is:
  \begin{definition}
  A graph $G$ is \emph{claw-free} if 
  $\overline{G[N(u)]}$ is triangle-free for 
  all $u\in V(G)$.
  \end{definition}
  \begin{lemma}\label{lemma:help0}  
  If $u$ is a vertex of a claw-free graph $G$,
  and $v\in N(u)$, then\break $G[N(v)\sm\nhs(u)]$
  is a complete graph.
  \end{lemma}
  \begin{proof}
  Let $w_1,w_2$ be two arbitrary 
  vertices of $G[N(v)\sm \nhs(u)]$.
  There are edges from $v$ to $w_1,w_2$ and $u$,
  and no edges from $u$ to $w_1$ and $w_2$. An
  edge between $w_1$ and $w_2$ is the only way
  to avoid a claw. 
  \end{proof}
\subsection{Independence complexes}
  \begin{definition}
  Let $G$ be a graph. The \emph{independence 
  complex
  of $G$}, $\ind(G)$ has vertex set $V(G)$ and
  its simplices are the independent subsets
  of $V(G)$.
  \end{definition}
  Some basic properties are:
  \begin{itemize}
  \bu If $U\subseteq V(G)$, then $\ind(G)[U]=
      \ind(G[U])$. 
  \bu If $u\in V(G)$ then $\ind(G\sm N(u))$ is
      a cone with apex $u$.
  \bu If $u\in V(G)$ and $\sigma\in\ind(G)$
      then there is a $v\in \nhs(u)$ such that
      $\sigma\cup \{v\}\in\ind(G)$.
  \bu If $u,v\in V(G)$ and $\{u,v\}$ is a 
      connected component of $G$, then
      $\ind(G)\simeq
      \susp(\ind(G\sm\{u,v\})).$
  \end{itemize}
  Two results from \cite{E} are needed. The
  proofs are short, so they are included for
  completeness.
  \begin{lemma}\label{lemma:fold}
  If $N(v)\subseteq N(w)$ then $\ind(G)$ 
  collapses onto $\ind(G\sm\{w\})$.
  \end{lemma}
  \begin{proof}
  Let $\{\sigma_1,\sigma_2,\ldots \sigma_k\}
  =\{\sigma\in\ind(G)\mid w\in\sigma, v\not\in
  \sigma\}$ be ordered such that if $\sigma_i
  \supseteq \sigma_j$ then $i<j$. The successive
  removals from $\ind(G)$
  of $\{\sigma_1,\sigma_1\cup\{v\}\},$
  $\{\sigma_2,\sigma_2\cup\{v\}\}, \ldots,
  \{\sigma_k,\sigma_k\cup\{v\}\}$ are elementary
  collapse steps.
  \end{proof}
  \begin{lemma}\label{lemma:rmCompN}
  If $u\in V(G)$ and $G[N(u)]$ is a complete
  graph, then
  \[\ind(G)\simeq\bigvee_
  {\hideS{v\in N(u)}}
  \susp(\ind(G\sm\nhs(v)))\]
  \end{lemma}
  \begin{proof}
  Let $\Delta_v=\ind(G\sm N(v))$ for all 
  $v\in \nhs(u)$. All $\Delta_v$ are 
  contractible, and $\Delta_{v_1}\cap
  \Delta_{v_2}\subseteq \Delta_u$ for all
  distinct $v_1,v_2\in N(u)$. By 
  Lemma~\ref{lemma:cf2}, and the
  third basic property of independence 
  complexes listed above,
  \[\ind(G)
  =\bigcup
  _{\hideS{v\in\nhs(u)}}
  %_{\makebox[0cm]{\scriptsize $v\in \nhs(u)$}} 
  \ind(G\sm N(v))
  =\bigcup
  _{\hideS{v\in \nhs(u)}}
  %_{\makebox[0cm]{\scriptsize $v\in \nhs(u)$}} 
  \Delta_v
  \simeq\bigvee
  _{\hideS{v\in N(u)}}
  %_{\makebox[0cm]{\scriptsize $v\in N(u)$}}
  \susp(\Delta_u\cap\Delta_v)
  =\bigvee
  _{\hideS{v\in N(u)}}
  %_{\makebox[0cm]{\scriptsize $v\in N(u)$}}  
  \susp(\ind(G\sm\nhs(v)))\]
  \end{proof}

\subsection{Higher connectivity}

\begin{figure}
  \begin{center}
  \includegraphics*{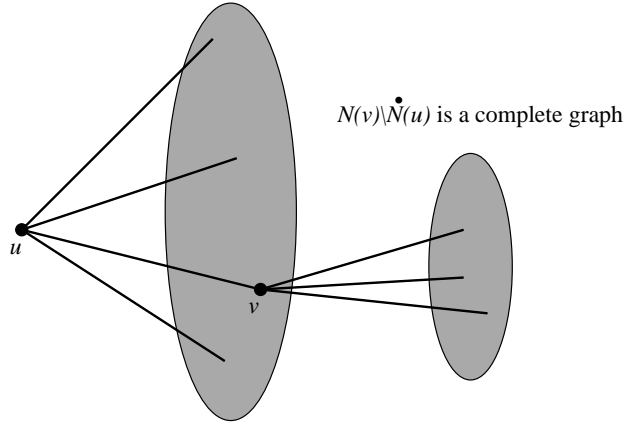}
  \end{center}
  \caption{The local structure of a claw-free graph}
  \label{fig:fig1}
\end{figure}

  Lemma~\ref{lemma:rmCompN} is a good tool
  for calculating the homotopy type of
  independence complexes of graphs where
  neighborhoods which form complete subgraphs
  can be found. In general this is not the
  case for claw-free graphs, but as illustrated
  in Figure~\ref{fig:fig1}, the
  situation is quite similar.
  It is probably impossible to use the
  local structure of claw-free graphs to
  calculate the homotopy type of their
  independence complexes recursively
  without running into devastating  
  identifications on the resulting
  topological space. However, in 
  Theorem~\ref{theorem:clawFreeMain} we
  show that the connectivity of independence
  complexes of claw-free graphs can be 
  handled.

  \begin{lemma}\label{lemma:help1}
  If $u,v\in V(G)$, $N(u)=\{v\}$, and
  $\ind(G\sm \nhs(v))$ is
  $(n-1)$--connected, then $\ind(G)$ is
  $n$--connected.
  \end{lemma}
  \begin{proof}
  The neighborhood of every vertex in 
  $N(v)\sm \{u\}$ contains $v$, and $v$ is the
  only vertex adjacent to $u$. Hence
  $\ind(G)$ collapses onto 
  $\ind(G\sm(N(v)\sm\{u\}))$ by repeated use of
  Lemma~\ref{lemma:fold}. The vertices $u$ and
  $v$ form a connected component of 
  $G\sm(N(v)\sm\{u\})$, so 
  $\ind(G\sm(N(v)\sm\{u\})) \simeq
   \susp(\ind((G\sm(N(v)\sm\{u\}))\sm \{u,v\}))$
  $=\susp(\ind(G\sm\nhs(v)))$. Since
  $\ind(G\sm\nhs(v))$ is $(n-1)$--connected,
  $\ind(G)$ is \mbox{$n$--connected.}
  \end{proof}

  \begin{lemma}\label{lemma:help2}
  Let $G$ be a graph with three vertices $u,v_1,$
  and $v_2$, such that\break 
  $\{v_1,v_2\}\not\in E(G)$,
  $N(u)=\{v_1,v_2\}$, and both 
  $G[N(v_1)\sm\{u\}]$ and $G[N(v_2)\sm\{u\}]$ are
  complete graphs. If $\ind(G\sm(\nhs(u)\cup 
  (N(v_1)\cap N(v_2))))$ is $(n-1)$--connected 
  and 
  $\ind(G\sm(\nhs(w_1)\cup\nhs(w_2)\cup\{u\}))$
  is $(n-2)$--connected for every
  $\{w_1,w_2\}\in
  E(\overline{G[N(v_1)\cup N(v_2)\sm\{u\}]}),$
  then $\ind(G)$ is $n$--connected.
  \end{lemma}

  \begin{proof}
  Let $H=G\sm (N(v_1)\cap N(v_2))$. First we 
  prove that $\ind(H)$ is $n$--connected, and 
  then the rest follows easily. 
  If $N_G(v_1)\subseteq N_G(v_2)$ then $v_1$
  is isolated in $H$ and $\ind(H)$ is a cone
  with apex $v_1$ and $n$--connected.
  Now assume that $N_G(v_1)\not\subseteq 
  N_G(v_2)$. The vertices
  $v_1$ and $v_2$ of $H$ have disjoint and
  complete neighborhoods, which fits good with
  using Lemma~\ref{lemma:rmCompN} twice,
  \[\ind(H)\simeq \bigvee
  %_{w_1\in N_H(v_1)}
  _{\hideS{w_1\in N_H(v_1)}}
  \susp(\ind(H\sm \nhs_H(w_1)))\]
  and 
  \[\ind(H\sm \nhs_H(w_1))=
  \bigvee
  %_{w_2\in N_H(v_2)\sm \nhs_H(w_1)}
  _{\hideS{w_2\in N_H(v_2)\sm \nhs_H(w_1)}}
  \susp(\ind(H\sm (\nhs_H(w_1)\cup\nhs_H(w_2)))).
  \]
  There is an edge between $w_1$ and $w_2$ in
  $\overline{G[N_G(v_1)\cup N_G(v_2)\sm\{u\}]}$ 
  if and only if 
  $w_1\in N_H(v_1)=N_G(v_1)\sm
  \{u\}$ and $w_2\in N_H(v_2)\sm \nhs_H(w_1)=
  (N_G(v_2)\sm\{u\})\sm \nhs_H(w_1)$.
  We assumed that  $\ind(H\sm (\nhs_H(w_1)\cup
  \nhs_H(w_2)))=\ind(G\sm (\nhs_G(w_1)\cup
  \nhs_G(w_2)\cup\{u\}))$ is $(n-2)$--connected 
  for every $\{w_1,w_2\}\in
  E(\overline{G[N(v_1)\cup N(v_2)\sm\{u\}]}),$
  therefore $\ind(H\sm \nhs_H(w_1))$ is 
  $(n-1)$--connected for every $w_1\in N_H(v_1)$.
  Well, actually not for all $w_1\in N_H(v_1)$
  because of that. If $\nhs_H(w_1)\supset 
  N_H(v_2)$, then we cannot use 
  Lemma~\ref{lemma:rmCompN} a second time, but
  then $\ind(H\sm \nhs_H(w_1))$ is a cone with
  apex $v_2$ and $(n-1)$--connected.

  All $\ind(H\sm \nhs_H(w_1))$ are 
  $(n-1)$--connected, so 
  $\ind(H)=\ind(G\sm(N_G(v_1)\cap N_G(v_2)))$ is 
  $n$-connected. The intersection of
  $\ind(G\sm(N_G(v_1)\cap N_G(v_2)))$ and
  $\ind(G\sm N_G(u))$ is $\ind(G\sm(\nhs_G(u)
  \cup (N_G(v_1)\cap N_G(v_2))))$ which is 
  assumed to
  be $(n-1)$-connected. $\ind(G\sm N_G(u))$ is
  a cone with apex $u$ and $n$--connected. Thus
  the union of $\ind(G\sm(N_G(v_1)\cap 
  N_G(v_2)))$ and $\ind(G\sm N_G(u))$, 
  $\ind(G\sm(N_G(v_1)\cap N_G(v_2) \sm \{u\}))$, 
  is $n$--connected. Finally, by repeated use
  of Lemma~\ref{lemma:fold}, $\ind(G)$ collapses
  onto 
  $\ind(G\sm(N_G(v_1)\cap N_G(v_2) \sm \{u\}))$
  since $N_G(w)\supset N_G(u)$ for all
  $w\in N_G(v_1)\cap N_G(v_2) \sm \{u\}$, and
  hence $\ind(G)$ is $n$--connected.
  \end{proof}

  \begin{theorem}\label{theorem:clawFreeMain}
  Let $u$ be a vertex of a claw-free graph $G$.
  If
  \begin{itemize}
  \bu $\ind(G\sm \nhs(v))$ is $(n-1)$--connected
  for every $v\in N(u)$ such that 
  \mbox{$\nhs(v)\supseteq \nhs(u)$},
%
%  WHY IS THIS ONE HERE?
%
%  \bu $\ind(G\sm (\nhs(v)\cup \nhs(u)))$ is 
%  $(n-2)$--connected
%  for every $v\in N(u)$ such that $\nhs(v)\
%  \not \supseteq \nhs(u)$,
  \bu $\ind(G\sm(\nhs(u)\cup(N(v_1)
  \cap N(v_2))))$ is $(n-1)$--connected
  for every $\{v_1,v_2\}\in 
  E(\overline{G[N(u)]})$,
  \bu $\ind(G\sm(\nhs(u)\cup\nhs(w_1)
  \cup\nhs(w_2)))$ is $(n-2)$--connected for
  every $\{w_1,w_2\}\in
  E(\overline{G[N(v_1)\cup N(v_2)\sm \nhs(u)]})$
  where $\{v_1,v_2\}\in E(\overline{G[N(u)]})$,
  \end{itemize}
  then $\ind(G)$ is $n$--connected.
  \end{theorem}
  \begin{proof}
  Define $\Delta_v=\ind(G\sm (N(u)\sm \{v\}))$ 
  for all $v\in N(u)$, and $\Delta_{v_1,v_2}=
  \ind(G\sm (N(u)\sm\{v_1,v_2\}))$ for all
  $\{v_1,v_2\} \in E(\overline{ G[N(u)] })$.

  Any face of $\ind(G)$ either contains a vertex
  from $\nhs(u)$ or can be extend with it.
  There is a face of $\ind(G)$ with two distinct 
  vertices $v_1,v_2$ of $\nhs(u)$ exactly when
  $\{v_1,v_2\}\in E(\overline{ G[N(u)] })$.
  But there can never be three vertices since
  the complement of a neighborhood in a 
  claw-free graph is triangle-free. A vertex
  $v$ of $N(u)$ such that $\nhs(v)\supseteq
  \nhs(u)$ can never be together with another
  vertex from $N(u)$ in a face of $\ind(G)$.
  We can cover $\ind(G)$:
  \[\ind(G)\wideEqual\bigcup_{\hideS{v\in N(u)
  \atop\nhs(v)\supseteq\nhs(u)}} \Delta_v\wideCup
  \bigcup _{\hideT{\{v_1,v_2\}\in 
  E(\overline{G[N(u)]})}}\Delta_{v_1,v_2}\]
  We will now show that the subcomplexes we
  cover with are $n$--connected and that their
  intersections are $(n-1)$--connected. From
  that we can conclude that $\ind(G)$ is
  $n$--connected by Lemma~\ref{lemma:cf1}. 
  The cases are:
  \begin{itemize}
  \item[(a)] $\Delta_v$ is $n$--connected for 
  all $v\in N(u)$ such that $\nhs(v)\supseteq
  \nhs(u)$.
  \item[(b)] $\Delta_{v_1,v_2}$ is $n$--connected
  for all $\{v_1,v_2\}\in E(\overline{G[N(u)]})$.
  \item[(c)] The intersection of at least two 
  different subcomplexes
  from (a) and (b) is $(n-1)$--connected:
  \begin{itemize}
  \item[(i)] One of the subcomplexes is a 
  $\Delta_v$.
  \item[(ii)] None of the subcomplexes is a
  $\Delta_v$, and there are two subcomplexes 
  $\Delta_{v_1,v_2}$ and $\Delta_{v_3,v_4}$ such
  that $\{v_1,v_2\} \cap  \{v_3,v_4\}=\emptyset$.
  \item[(iii)] The subcomplexes are
  $\Delta_{v,v_1},\Delta_{v,v_2},\ldots 
  \Delta_{v,v_k}.$
  \end{itemize}
  \end{itemize}

  \textbf{Case a.}
  Let $v$ be a vertex of $N(u)$ such that
  $\nhs(v)\supseteq\nhs(u)$. The neighborhood
  of $u$ in $G\sm(N(u)\sm\{v\})$ is $\{v\}$, so
  by Lemma~\ref{lemma:help1}, 
  $\Delta_v=\ind(G\sm(N(u)\sm\{v\}))$ is
  $n$--connected since
  \[\ind( (G\sm(N(u)\sm\{v\}))\sm 
  \nhs_{G\sm(N(u)\sm\{v\})}(v))
  =\ind(G\sm\nhs(v))\]
  is $(n-1)$--connected by assumption.

  \textbf{Case b.}
  Let $\{v_1,v_2\}$ be an edge of 
  $\overline{G[N(u)]}$ and define 
  $H=G\sm(N(u)\sm\{v_1,v_2\})$. We are to prove
  that $\Delta_{v_1,v_2}=\ind(H)$ is 
  $n$-connected, and we will us 
  Lemma~\ref{lemma:help2} to do that. Let's
  check the conditions of the lemma. The
  three vertices we use are $u,v_1,v_2$.
  \begin{itemize}
  \bu $\{v_1,v_2\}\not\in E(H)$.
  \bu $N_H(u)=\{v_1,v_2\}$.
  \bu By Lemma~\ref{lemma:help0}, $H[N_H(v_1)
  \sm \{u\}]=G[N_G(v_1)\sm\nhs_G(u)]$ is
  a complete graph.
  \bu By the same reason $H[N_H(v_2)\sm \{u\}]$
  is a complete graph.
  \bu From the inclusions $H\subset G$ and
  $N_G(u)\sm \{v_1,v_2\} \subseteq \nhs_G(u)
  \cup(N_G(v_1)\cap
  N_G(v_2))) $ we get that
  $\ind(H \sm (\nhs_H(u)\cup(N_H(v_1)\cap
  N_H(v_2)))) = 
  \ind(G \sm (\nhs_G(u)\cup(N_G(v_1)\cap
  N_G(v_2))))$ which is $(n-1)$--connected
  by assumption.
  \bu In the same way 
  $\ind(H\sm(\nhs_H(w_1)\cup\nhs_H(w_2)
  \cup\{u\}))$
  is $(n-2)$--connected for every
  $\{w_1,w_2\}\in
  E(\overline{H[N_H(v_1)\cup N_H(v_2)\sm\{u\}]}),$
  since $\ind(G\sm(\nhs_G(u)\cup\nhs_G(w_1)
  \cup\nhs_G(w_2)))$ is $(n-2)$--connected for
  every\\ $\{w_1,w_2\}\in
  E(\overline{G[N_G(v_1)\cup N_G(v_2)\sm 
  \nhs_G(u)]})$
  where $\{v_1,v_2\}\in E(\overline{G[N_G(u)]})$
  by assumption.
  \end{itemize}

  \textbf{Case c.}
  First note that the intersection with any of
  the subcomplexes $\Delta_v$ and 
  $\Delta_{v_1,v_2}$ with $\ind(G\sm N(u))$ is 
  $\ind(G\sm N(u))$. And that is a cone with
  apex $u$ and thus contractible. After
  sufficient many intersections of subcomplexes
  we will see that one ends up with $\ind(G\sm 
  N(u))$
  for which the connectedness is allright.

  \textbf{Case c.i.}
  Say that one of the subcomplexes is 
  $\Delta_{v_1}$. If $v_1\neq v_2$ then 
  $\Delta_{v_1}\cap\Delta_{v_2}= \ind(G\sm N(u))$.
  If $\{v_2,v_3\}\in E(\overline{G[N(u)]})$
  and $\nhs(v_1)\supseteq \nhs(u)$ then
  $v_1\not\in \{v_2,v_3\}$ and
  $\Delta_{v_1}\cap\Delta_{v_2,v_3}= 
  \ind(G\sm N(u))$. We conclude that
  a intersection where one of the subcomplexes
  is $\Delta_{v_1}$ is $(n-1)$--connected.

  \textbf{Case c.ii.}
  The intersection of two subcomplexes
  $\Delta_{v_1,v_2}$ and $\Delta_{v_3,v_4}$ such that
  $\{v_1,v_2\} \cap  \{v_3,v_4\}=\emptyset$
  is $\ind(G\sm N(u))$ so the complete
  complete intersection is also $\ind(G\sm N(u))$
  which is $(n-1)$--connected.

  \textbf{Case c.iii.}
  $\cap_{i=1}^k \Delta_{v,v_i}=\ind(G\sm (N(u)\sm v))$.
  We assumed that $\ind(G\sm (\nhs(v)\cup \nhs(u)))$ 
  is $(n-2)$--connected
  for every $v\in N(u)$ such that $\nhs(v)\
  \not \supseteq \nhs(u)$. By Lemma~\ref{lemma:help1},
  $\ind(G\sm (N(u)\sm v))$ is $(n-1)$--connected.
  \end{proof}

\section{Asymptotic higher connectivity}
  It was proved in \cite[Theorem 26]{E}
  that for any graph
  $G$ with maximal degree $d$, $\ind(G)$ is
  $(\lfloor (n-2d-1)/2d \rfloor)$--connected,
  where $d$ is the maximal degree of a vertex
  of $G$. For a graph property, it is an
  interesting task to find the best $c$, such
  that for $G$ with the property,
  $\ind(G)$ is $f(n,d)$--connected where
  $f(d,d)$ grows asymptotically as $cn/d$. In
  \cite{E,ST} it was proved that $c=1/2$ if we
  put no restriction on the graphs. In this
  section we prove that $c\geq 2/3$ for
  claw-free graphs.

\begin{lemma}\label{lemma:clawFreeIneq}
  If $G$ is a claw-free graph
  with maximal degree $d$, $u\in V(G)$, and
  $\{v_1,v_2\} \subseteq N(u)$ but
  $\{v_1,v_2\} \not \in E(G)$, then
  \[ \# \nhs(u)\cup(N(v_1)\cap N(v_2)) \leq 
  \lfloor (3d+2)/2 \rfloor \]
\end{lemma}
\begin{proof}
  For every vertex in the 
  neigborhood of $u$ other than $v_1$ and 
  $v_2$, at least one of $v_1$ and $v_2$ must have 
  an edge to it since $G$ is claw-free. Therefore
  either $v_1$ or $v_2$ must have edges
  to at least half of the elements of
  $N(u)\sm \{v_1,v_2\}.$ Assume that it
  is $v_1$. Insert
  \[ \# N(u) \cap N(v_1) \geq \lceil 
  (\# N(u) -2)/2 \rceil 
  \Rightarrow
   \# \nhs(u) \cap N(v_1) \geq \lceil 
   \# N(u)/2 \rceil \]
   into
   \[\begin{array}{rcl} 
     \# N(v_1)\cap N(v_2) \sm \nhs(u) &
     \leq & \# N(v_1)\sm\nhs(u) \\
     & = & \# N(v_1) - \# N(v_1) \cap \nhs(u) \\
     & \leq & \# N(v_1) - \lceil \# N(u)/2 \rceil
   \end{array}\]
   to conclude that
   \[\begin{array}{rcl}
     \# \nhs(u)\cup(N(v_1)\cap N(v_2)) & = &
     \# \nhs(u) + \# N(v_1)\cap N(v_2)\sm \nhs(u) \\
     & \leq & \# \nhs(u) + \# N(v_1) - \lceil \# N(u)/2 \rceil \\
     & = & 1 + \# N(u) + \# N(v_1) - \lceil \# N(u)/2 \rceil \\
     & = & 1 + \# N(v_1) + \lfloor \# N(u)/2 \rfloor \\
     & \leq & 1 + d + \lfloor d/2 \rfloor \\
     & = & \lfloor (3d+2)/2 \rfloor \\
   \end{array}\]
\end{proof}

\begin{theorem}
  If $G$ is a claw-free graph with $n$
  vertices and maximal degree $d$, then
  $\ind(G)$ is
  $\lfloor (2n-1)/(3d+2)-1 \rfloor$--connected.
\end{theorem}
\begin{proof}
  If $d=0$ the statement is true, so assume that $d\geq 1$.
  If $0<n\leq (3d+2)/2$ the statement is that $\ind(G)$
  is $(-1)$--connected. This means that the complex
  is nonempty, which is true. The proof is by induction
  over the number of vertices. Note that subgraphs of $G$
  never have higher maximal degree than $d$.

  Assume that $n>(3d+2)/2$ and fix a vertex $u$ of $G$.
  The independence complex of $G$ is broken up into
  smaller pieces with bounded connectivity and patched
  together with Theorem~\ref{theorem:clawFreeMain}.
  The next step is to check that the conditions of
  the theorem are fullfilled.
  \begin{itemize}
  \bu Let $v$ be a vertex in $N(u)$. There are at most
  $d+1$ elements in $\nhs(v)$, and $(3d+1)/2 \geq d+1$,
  so $\ind(G\sm \nhs(v))$ is  
  $(\lfloor (2n-1)/(3d+2)-1 \rfloor-1)$--connected by
  induction.
  \bu By Lemma~\ref{lemma:clawFreeIneq} 
  $\# \nhs(u)\cup(N(v_1)\cap N(v_2)) \leq 
  \lfloor (3d+2)/2 \rfloor$ for every $\{v_1,v_2\}$ in
  $E(\overline{G[N(u)]})$. Thus $\ind(G\sm( \nhs(u)
  \cup(N(v_1)\cap N(v_2))))$ is $(\lfloor (2n-1)/
  (3d+2)-1 \rfloor-1)$--connected by induction.
  \bu For every $\{w_1,w_2\}\in
  E(\overline{G[N(v_1)\cup N(v_2)\sm \nhs(u)]})$
  where $\{v_1,v_2\}\in E(\overline{G[N(u)]})$,
  the intersection of $\nhs(u)$ and
  $\nhs(w_1)\cup\nhs(w_2)$ contains $v_1$ and $v_2$,
  so $\# \nhs(u)\cup\nhs(w_1)\cup\nhs(w_2) \leq 3d+1$.
  Therefore $\ind(G\sm(\nhs(u)\cup\nhs(w_1)
  \cup\nhs(w_2)))$ is $(\lfloor (2n-1)/
  (3d+2)-1 \rfloor-2)$--connected by induction.
  \end{itemize}
  We conclude by Theorem~\ref{theorem:clawFreeMain}
  that $\ind(G)$ is
  $\lfloor (2n-1)/(3d+2)-1 \rfloor$--connected.
\end{proof}

\section{Connectivity of $\mathcal{C}^k_n$}
  We will treat two classes of independence
  complexes of claw-free graphs
  introduced by Kozlov~\cite{K}.
  Let $L_n^k$ be the graph with vertex set
  $\{1,2,\ldots, n\}$ and two vertices $i<j$
  are adjacent if $j-i<k$. Define
  $\mathcal{L}^k_n=\ind(L_n^k)$.
  For $n\leq 0$ let $\mathcal{L}^k_n=\emptyset$.
  In Engstr\"om~\cite[Corollary 21]{E}
  it was proved that
  \[\mathcal{L}_n^k
  \simeq \bigvee_{1\leq i< \min \{k,n\}}
  \susp(\mathcal{L}_{n-k-i}^k)\]
  using something like Lemma~\ref{lemma:rmCompN}.
  It follows directly that $\mathcal{L}_n^k$
  is $l_{n,k}$--connected, where
  \[l_{n,k}=\left\lfloor \frac{n-1}{2k-1}
   -1 \right\rfloor.\]

  The second class is build from $C_n^k$ which
  is a graph with vertex set $\{1,2,\ldots, n\}$ 
  and two vertices $i<j$ are adjacent if $j-i<k$
  or $(n+i)-j<k$. Define $\mathcal{C}_n^k
  =\ind(C_n^k)$. The homotopy type of $\mathcal
  {C}_n^2$ was determined in \cite{K}, and used
  by Babson and Kozlov in their proof of 
  Lov\'asz conjecture \cite{BK}. Some other cases
  where treated in \cite{E}, but in general the
  homotopy type of $\mathcal{C}_n^k$ is not 
  known. Removing at least $k$ consecutive 
  vertices from $\mathcal{C}_n^k$ gives a complex
  of the $\mathcal{L}$ type which we know
  the higher connectivity of. We will cover
  $\mathcal{C}_n^k$ with $\mathcal{L}$ type
  complexes and then use 
  Theorem~\ref{theorem:clawFreeMain} to bound 
  the connectivity of it. Why is $C_n^k$
  claw-free? If we for example pick three
  elements of $N(k)$, then two of them must
  be either larger or smaller than $k$,
  which forces their difference
  smaller than $k$, and they are adjacent.
\begin{theorem}
  If $n\geq 6(k-1)$ then $\mathcal{C}_n^k$ is
  $c_{n,k}$--connected, where
  \[c_{n,k}=\left\lfloor \frac{n+1}{2k-1} -2
  \right\rfloor.\]
\end{theorem}
\begin{proof}
  We are to check the conditions of 
  Theorem ~\ref{theorem:clawFreeMain}. Let
  $u=3k-2$.
\begin{itemize}
  \bu There no $v\in N(u)$ such that 
  $\nhs(v)\subseteq \nhs(u)$.

  \bu If $\{v_1,v_2\}\in E(\overline{
  G[N(u)]})$ then $N(v_1)\cap N(v_2)\subseteq
  \nhs(u)$, so  $\ind(G\sm(\nhs(u)\cup(N(v_1)
  \cap N(v_2))))= \ind(G\sm(\nhs(u)))
  \simeq \mathcal{L}_{n-(2k-1)}^k$ which is
  $l_{n-(2k-1),k}$--connected. Clearly
  $c_{n,k}-1 \leq l_{n-(2k-1),k}$.

  \bu Choose $v_1=2k-1, v_2=4k-3,
  w_1=k,$ and $w_2=5k-4$ to minimize the
  size of $\ind(C_n^k\sm (\nhs(u)\cup
  \nhs(w_1)\cup\nhs(w_2)))\simeq 
  \mathcal{L}_{n-(6k-5)}^k$ which is 
  $l_{n-(6k-5),k}$--connected. Clearly
  $c_{n,k}-2 = l_{n-(6k-5),k}$.
\end{itemize} 
\end{proof}

\bibliographystyle{amsplain}

\begin{thebibliography}{10}
  \bibitem{BK} E. Babson, D.N. Kozlov, \textit{
  Proof of
  the Lov\'asz conjecture}, Ann. of Math. (2),
  to appear.  

  \bibitem{B} A. Bj\"orner, 
  \textit{Topological Methods,} 
  in: ``Handbook of Combinatorics'' (eds. R. 
  Graham, M. Gr\"otschel, and L. Lov\'asz), 
  North-Holland, 1995, 1819--1872.

  \bibitem{BLVZ} A. Bj\"orner, L. Lov\'asz, S.T. 
  Vre\'cica, R.T. \v{Z}ivaljevi\'c, 
  \textit{Chessboard complexes and matching 
  complexes}, J. London Math. Soc. (2) 
  \textbf{49} (1994), 25--49.
   
  \bibitem{E} A. Engstr\"om,
  \textit{Complexes of Directed Trees and
  Independence Complexes},
  \verb|http://www.arxiv.org/abs/math/0508148|

  \bibitem{K} D.N. Kozlov, \textit{Complexes of
  directed trees}, J. Combinat. Theory Ser. A 88
  (1999), no. 1, 112-122.

  \bibitem{ST} T. Szab\'o, G. Tardos,
  \textit{Extremal problems for transversals in 
  graphs with bounded degree}, Combinatorica, 
  to appear.

\end{thebibliography}

\end{document}